\documentclass[12pt]{article}
\usepackage{amsmath,amssymb, amsthm, latexsym, amsfonts}
\usepackage{graphicx}
\usepackage{subcaption}

\newtheorem{thm}{Theorem}[section]

\newtheorem{lemma}{Lemma}[section]

\usepackage{multicol}
\usepackage{vmargin, marginnote}
\usepackage[mathscr]{eucal}
\usepackage{xcolor}

\newcommand{\beq}{\begin{equation*}}

\newcommand{\eeq}{\end{equation*}}

\newcommand{\Beq}{\begin{equation}}

\newcommand{\Eeq}{\end{equation}}

\title{Existence and numerical approximation of a one-dimensional Boussinesq system with variable coefficients on a finite interval}

\author{Deissy Marcela Pizo\footnote{Universidad del Valle, Departamento de Matem\'aticas, Cali, Colombia, email: deissy.pizo@correounivalle.edu.co }, Juan Carlos Mu\~noz Grajales \footnote{Universidad del Valle, Departamento de Matem\'aticas, Cali, Colombia, Corresponding author, email: juan.munoz@correounivalle.edu.co. } }

\begin{document}

\maketitle



\begin{abstract}
In this paper, we investigate the well-posedness of a nonlinear dispersive model with variable coefficients that describes the evolution of surface waves propagating through a one-dimensional shallow water channel of finite length with irregular bottom topography. To complement the theoretical analysis, we utilize the numerical solver developed by the authors in \cite{PizoMunoz} to approximate solutions of the model on a finite spatial interval, considering various parameter values and forms of the variable coefficients in the Boussinesq system under study. Additionally, we present preliminary numerical experiments addressing an inverse problem: the reconstruction of the initial wave elevation and fluid velocity from measurements taken at a final time. This is achieved by formulating an optimization problem in which the initial conditions are estimated as minimizers of a functional that quantifies the discrepancy between the observed final state and the numerical solution evolved from a trial initial state.

\end{abstract}

\vspace{0.5cm}



\section{Introduction}

Boussinesq systems have been extensively studied in recent years from both analytical and numerical perspectives, as models for water wave propagation in shallow channels with constant or variable depth \cite{Bona1998, Bona2002, Bona2004, Grajales2004}. In particular, issues such as well-posedness, existence, and the orbital stability of traveling wave solutions have been addressed using a variety of mathematical tools, including functional analysis, Sobolev space theory, variational methods, and numerical techniques.

In this paper, we begin by investigating the well-posedness of a Boussinesq-type system with variable coefficients, given by
\begin{equation}\label{Boussi_intro}
	\begin{cases}
		\begin{array}{c l}
			N_t+\left[ (1+\alpha c^{2}\left( \xi \right)N)   V\right] _{\xi}-\frac{\beta}{6}N_{\xi\xi t}&=0 \:\:\:,\:\:\left(\xi,t\right)\in \left[0,L\right]\times\left[0,T\right].\\
			V_t+(c(\xi) N)_{\xi}+\dfrac{1}{2}\alpha (c^{2}V^{2})_{\xi }-\frac{\beta}{6}V_{\xi\xi t}&=0,
		\end{array}
	\end{cases}
\end{equation}
subject to the initial conditions
\begin{equation}\label{Inicial_intro}
	\begin{array}{c l}
		N(\xi,0)&=N_0(\xi),\\
		V(\xi,0)&=V_0(\xi),\:\:\: \xi \in \left[0,L\right],
	\end{array}
\end{equation}
and Dirichlet-type boundary conditions
\begin{equation}\label{Frontera_intro}
	\begin{array}{c l}
		N(0,t)=N(L,t)=0,\:\: t\in \left[0,T\right]\\
		V(0,t)=V(L,t)=0,\:\: t\in \left[0,T\right].
	\end{array}
\end{equation}

This system models the evolution of surface waves propagating through a one-dimensional shallow water channel of finite length with irregular bottom topography.

To complement the theoretical analysis, we make use of the numerical solver developed by the authors in \cite{PizoMunoz} to approximate solutions of the system \eqref{Boussi_intro} on a bounded spatial domain, considering various parameter values and functional forms for the coefficient $c(\xi)$. Furthermore, we present preliminary numerical experiments that address an inverse problem: reconstructing the initial wave elevation and fluid velocity from measurements taken at a final time. This problem is formulated as an optimization task, where the initial conditions are identified as minimizers of a functional that quantifies the discrepancy between the observed final state and the numerically simulated state evolved from a trial initial condition.

To the best of our knowledge, this type of inverse problem for a dispersive system such as \eqref{Boussi_intro} has not previously been studied and therefore constitutes the main contribution of this work.

The remainder of the manuscript is organized as follows. Section 2 introduces the notation and theoretical tools that will be used throughout the paper. Section 3 is devoted to establishing the existence, uniqueness, and continuous dependence on initial data for solutions of the initial-boundary value problem \eqref{Boussi_intro}–\eqref{Frontera_intro}, using Green’s function techniques associated with the spatial second-order spatial operator. Section 4 presents the numerical framework and results regarding the inverse problem described above, using the solver introduced in \cite{PizoMunoz}. Section 5 concludes the paper with a discussion of the main findings and potential directions for future work.

\section{Preliminaries}

In this chapter, we recall some notation, definitions, and standard results from Functional Analysis that will be used throughout the manuscript.

Let $H$ be a Hilbert space. We denote the inner product and the norm in $H$ by $\langle \cdot , \cdot \rangle_H$ and $|\cdot|_H$, respectively. 



The Hilbert space of all measurable functions 
$v:\left[0,L\right] \rightarrow \mathbb{R}$ that are square-integrable over the interval $\left[0, L \right]$ is denoted by $L^2$, that is, $L^2=L^2\!\left(\left[0,L \right]  \right)$. This space is endowed with the standard inner product and the asscociated norm
\begin{equation}\label{Ana6}
\left\langle w,v\right\rangle_{L^2}=\int_{0}^{L}w(x)v(x)\,dx, 
\qquad 
\left\|v \right\|_{L^2} =\left(\int_{0}^{L}\left( v(x)\right)^2 dx\right)^ {\frac{1}{2}}.
\end{equation}
We denote by $H^1$ the Sobolev space consisting of all measurable functions $f:\left[0,L\right] \rightarrow \mathbb{R}$ whose weak derivative $f'$ existes and satisfies $f, f' \in L^2(0,L)$. In this work, the space $H^1$ is equipped with the inner product
\begin{equation}
\left\langle f,g\right\rangle_{H^1}=\int_{0}^{L}f(x) g(x)\, dx + \frac{\beta}{6}\int_{0}^{L}f^{'}(x)g^{'}(x)\,dx,
\end{equation}
and the corresponding norm
\begin{equation}\label{MA27}
\left\|f \right\|_{H^1} =\left(\int_{0}^{L}\left(f(x)\right)^2 dx + \frac{\beta}{6}\int_{0}^{L}\left( f^{'}(x)\right)^2 dx\right)^{\! \frac{1}{2}}.
\end{equation}

We recall below some fundamental results concerning the function spaces introduced above, as they will be used frequently in the subsequent chapters.

For $f \in  H^1$, the inequality
\begin{equation}
\left\|f^{'} \right\|_{L^2}\leq \frac{1}{\left(\frac{\beta}{6} \right)^{1/2} }\left\|f \right\|_{H^1}
\end{equation}
holds.  

An important result ensuring the compactness of the embedding $H^1([0, L]) \hookrightarrow L^{2}([0, L])$ is the Rellich-Kondrachov Theorem, which we state below.

\begin{thm}
The embedding of the Sobolev space $H^1([0, L])$ into $L^{2}([0, L])$ is compact.
\end{thm}

A proof of this result can be found in Folland (\cite{Folland}), p.~225. 
\\
The subspace of functions $f \in H^1([0, L])$ such that $f(0) = f(L) = 0$ (in the sense of traces) is denoted by $H_{0}^{1}([0, L])$. The continuous embedding of $H_0^{1}([0, L])$ into $L^{\infty}([0, L])$ follows from the Fundamental Theorem of Calculus, since
\begin{equation}
\left|f(x)\right|=\left|\int_{0}^{x}f^{'}(s)\,ds\right| \leq \int_{0}^{L} \left|f^{'}(s) \right|\,ds.
\end{equation}
Therefore, applying the Cauchy--Schwarz inequality, one deduces
\begin{equation}
\qquad \left\|f \right\|_{L^{\infty}}\leq \frac{L^{1/2}}{\left(\tfrac{\beta}{6} \right)^{1/2} }\left\|f \right\|_{H^1}.
\end{equation}

If $X$ is a Banach space with norm $\|\cdot\|_{X}$, and $0<T\leq \infty$, the space $\mathcal{L}_T^p(X)= L^{p}(0,T;X)$, with $1\leq p <\infty$, consists of all measurable functions $f:\left[0,T \right]\rightarrow X$ such that the norm
\begin{equation}
\| f\|_{\mathcal{L}^{p}_{T}}=\left(\int_{0}^{T}\|f(t) \|^{p}_{X}\, dt\right)^{1/p}
\end{equation}
is finite. Similarly, the space $C([0,T];X)=\mathcal{C}_{T}(X)$ consists of all continuous functions $f:[0,T]\rightarrow X$ such that the norm
\begin{equation}
\| f\|_{\mathcal{C}_{T}(X)}=\max_{t\in [0,T]}\| f(t)\|_{X} 
\end{equation}
is finite.



In what follows, we adopt the notation
\begin{equation*}
	\mathcal{L}_{T}^2=L^2\!\left( 0,T; H_{0}^1([0,L]) \right),\quad
	\mathcal{C}_T= C\!\left( 0,T; H_{0}^1([0,L]) \right),\quad
	L_{T}^2=L^2\!\left( 0,T; L^2([0,L]) \right),
\end{equation*}
where $\mathcal{L}_{T}^2$ is a Banach space with the norm
\begin{equation}
	\|N \|_{\mathcal{L}_{T}^2}=\left( \int_{0}^{T}\|N(t) \|_{H_{0}^1([0,L])}^2\, dt \right)^{1/2}.
\end{equation}
Here, the space $H_{0}^1$ is considered with the norm given in~\eqref{Ana6}. Similarly, the spaces $\mathcal{C}_T$ and $L_{T}^2$ are Banach spaces with the norms
\begin{equation}
	\|N \|_{\mathcal{C}_{T}}=\max_{t\in [0,T]}\|N(t) \|_{H_{0}^1([0,L])}, 
	\qquad 
	\|N \|_{L_{T}^2}=\left( \int_{0}^{T}\|N(t) \|_{L^2([0,L])}^2\, dt \right)^{1/2},
\end{equation}
respectively.

\section{Existence and uniqueness}

Let $L$, $\beta$, $T$ be positive constants, and $c(\xi), \alpha(\xi) \in L^{2}([0,L])$. Let us consider the system
\begin{equation}\label{Dee}
	\begin{cases}
		\begin{array}{c l}
			N_t+\left[ (1+\alpha c^{2}\left( \xi \right)N)   V\right] _{\xi}-\frac{\beta}{6}N_{\xi\xi t}&=0 \:\:\:,\:\:\left(\xi,t\right)\in \left[0,L\right]\times\left[0,T\right].\\
			V_t+(c(\xi) N)_{\xi}+\dfrac{1}{2}\alpha (c^{2}V^{2})_{\xi }-\frac{\beta}{6}V_{\xi\xi t}&=0,
		\end{array}
	\end{cases}
\end{equation}
subject to the initial conditions
\begin{equation}\label{Maa}
	\begin{array}{c l}
		N(\xi,0)&=N_0(\xi),\\
		V(\xi,0)&=V_0(\xi),\:\:\: \xi \in \left[0,L\right],
	\end{array}
\end{equation}
and Dirichlet-type boundary conditions
\begin{equation}\label{BBMLin33}
	\begin{array}{c l}
		N(0,t)=N(L,t)=0,\:\: t\in \left[0,T\right]\\
		V(0,t)=V(L,t)=0,\:\: t\in \left[0,T\right].
	\end{array}
\end{equation}

The existence of solutions to the system \eqref{Dee} will be studied
 using the method of Green's functions. This approach allows us to reformulate the problem as one of finding a fixed point of a suitably defined nonlinear operator. The existence and uniqueness of such a fixed point will then follow from an application of the Banach Fixed Point Theorem. 
 
Observe that the equations in the system \eqref{Dee} can be rewritten in the form
\begin{equation}
	\begin{array}{c l}
		\left(I-\frac{\beta}{6}\partial_{\xi}^2 \right) N_t =  --\partial_\xi \left[ (1+\alpha c^{2}\left( \xi \right)N)   V\right]\\
		\left(I-\frac{\beta}{6}\partial_{\xi}^2 \right) V_t=-\partial_\xi\left[ \left( c(\xi) N\right)+\dfrac{1}{2}\alpha (c^{2}V^{2})\right] .
	\end{array}
\end{equation}
Let $G(\xi,s)$ denote the Green's function associated with the differential operator $P=I-\frac{\beta}{6}\partial_{\xi}^2 $, which is given by
\begin{equation}\label{FuncionG}
	G(\xi, s) = \frac{1}{2 \sqrt{\beta/6}} \frac{ \cosh\Big( \frac{1}{\sqrt{\beta/6}} (L - |s - \xi | ) \Big) - \cosh \Big( \frac{1}{\sqrt{\beta/6}} (L - (\xi + s) ) \Big)}{\sinh\Big( \frac{L}{\sqrt{\beta/6} } \Big)}.
\end{equation}

It can be shown that the derivative of $G(\xi, s)$ with respect to the second variable $s$ is given by

\begin{equation}\label{DerivadaG}
	K(\xi,s) :=  \frac{\partial G(\xi, s)}{\partial s} 
  = \frac{3}{\beta} \left( \frac{ \sinh \Big( \frac{L-\xi-s}{\sqrt{\beta/6}} \Big) + \text{sign} (\xi- s) \sinh \Big( \frac{L- |\xi-s|}{\sqrt{\beta/6}} \Big)}{ \sinh \Big( \frac{L}{\sqrt{\beta/6}} \Big) }  \right),\xi \neq s.
\end{equation}
Note that
\[
\lim_{s \to \xi^-} K(\xi, s) - \lim_{s \to \xi^+} K(\xi,s) = \frac{6}{\beta}.
\]
It follows that a solution $\left(N,V \right) $ of the system $\left(\ref{Dee}\right)$ must satisfy the system
\begin{eqnarray}\label{Ma11}
	N_t & = &-\int_{0}^{L}G(\xi,s)\left[ (1+ \alpha c^{2}N\left( s,t\right) ) V\left( s,t\right) \right]_s ds \\
	& = & \int_{0}^{L}G_{s}(\xi,s)(1+ \alpha\left(s\right) c^{2}N\left( s,t\right) ) V\left( s,t\right)ds, \label{Ma11new}
\end{eqnarray}

\begin{eqnarray}\label{Ma1}
	V_t & = &-\int_{0}^{L}G(\xi,s) \left[ ( c\left(s \right)N\left( s,t\right)) + \dfrac{1}{2}\alpha \left( c^{2}V^{2}\left( s,t\right)\right)\right] _s ds \\
	& = & \int_{0}^{L}G_s(\xi,s)   \left[ ( c\left(s \right)N\left( s,t\right)) + \dfrac{1}{2}\alpha  \left( c^{2}V^{2}\left( s,t\right)\right)\right] ds.  \label{Ma1new}
\end{eqnarray}

The following lemma establishes important properties concerning the functions $G(\xi, s)$ and $K(\xi,s)$ defined above that are necessary in what follows.
\

\begin{lemma}\label{MA51}
	Let $\Phi_1$ and $\Phi_2$ be the operators defined by
	\begin{equation*}
		\Phi_1(\phi)(\xi)=\int_{0}^{L}G(\xi,s)\phi(s)ds \quad y \quad \Phi_2(\phi)(\xi)=\int_{0}^{L}K(\xi,s)\phi(s)ds.
	\end{equation*}
	Then, $\Phi_1$ and $\Phi_2$ are linear and continuous operators of $L^{2}$ on $H^{1}$, i.e, there exist positive constants $C_1$ and $C_2$ such that
	\begin{equation}\label{MA50}
		\left\| \Phi_1(\phi)\right\|_{H^{1}}\leq C_1\left\| \phi\right\|_{L^{2}}  \quad y \quad \left\| \Phi_2(\phi)\right\|_{H^{1}}\leq C_2\left\| \phi\right\|_{L^{2}}, 
	\end{equation}
	for all $\phi \in L^{2}$.
\end{lemma}

\textbf{Proof:} 
The linearity of the operators $\Phi_1$, $\Phi_2$ follows directly from the linearity of the integral operator involved. On the other hand, by applying the Cauchy-Schwarz inequality, we obtain
\begin{equation}\label{MA47}
	\begin{split}
		\left\| \Phi_1(\phi)\right\|^{2}_{L^{2}}=\int_{0}^{L} \left\lbrace \int_{0}^{L} G(\xi,s)\phi(s)ds \right\rbrace ^{2} dx \leq \left\| G\right\| ^{2}_{L^{2}(0,L;L^{2}[0,L])}\left\|\phi \right\| ^{2}_{L^{2}}
	\end{split}
\end{equation}
Furthermore, note that
\begin{equation}
	\begin{split}
		( \Phi_1(\phi))^{'}(\xi)=\int_{0}^{L} \partial_{\xi} G(\xi,s)\phi(s)ds = \int_{0}^{L} K(\xi,s)\phi(s)ds,
	\end{split}
\end{equation}

from where it is obtained

\begin{equation}\label{MA48}
	\begin{split}
		\left\|   	( \Phi_1(\phi))^{'} \right\|^{2}_{L^{2}}=\left\| \Phi_2(\phi)\right\|^{2}_{L^{2}}=\int_{0}^{L} \left\lbrace \int_{0}^{L} K(\xi,s)\phi(s)ds \right\rbrace ^{2} dx 
        \\\\ \leq \left\| K\right\| ^{2}_{L^{2}(0,L;L^{2}[0,L])}\left\|\phi \right\| ^{2}_{L^{2}}
 	\end{split}       
\end{equation}
On the other hand, the function $K$ defined in (\ref{DerivadaG}) is continuous in $[0,L]\times [0,L]$, except at $\xi=s$, where there is a jump discontinuity of the form
\[
\lim_{s \to \xi^-} K(\xi, s) - \lim_{s \to \xi^+} K(\xi,s) =K(\xi,\xi^{-})-K(\xi,\xi^{+}) =\frac{6}{\beta}.
\]
With this in mind, a direct calculation allows us to establish that
\begin{equation}
	\begin{split}
		( \Phi_2(\phi))^{'}(\xi)=[K(\xi,\xi^{-})-K(\xi,\xi^{+}) ]\phi(\xi)+\int_{0}^{L} \partial_{\xi} K(\xi,s)\phi(s)ds
		\\\\= \frac{6}{\beta} \phi(\xi)+\int_{0}^{L} K(\xi,s)\phi(s)ds.
	\end{split}
\end{equation}
Consequently, by the Cauchy-Schwarz inequality, we have
\begin{equation}\label{MA49}
	\begin{split}
		\left\|   	( \Phi_2(\phi))^{'} \right\|^{2}_{L^{2}} \leq \frac{2.6}{\beta^{2}}\left\| \phi\right\| ^{2}_{L^{2}} + 2 \left\| \int_{0}^{L}\partial_{\xi}K(\cdot,s)\phi(s)ds\right\|^{2}_{L^{2}}  \\\\ \leq 
          2\left(\frac{6}{\beta^{2}}+ \left\| \partial_{\xi}K\right\| ^{2}_{L^{2}(0,L;L^{2}[0,L])}  \right) \left\| \phi\right\|^{2}_{L^{2}}. 
	\end{split}
\end{equation}
Therefore, inequalities (\ref{MA47}), (\ref{MA48}),(\ref{MA49}) imply the desired results in (\ref{MA50}).
$\hfill\square$

The following theorem establishes a local well-posedness result for the general initial-boundary value problem associated with the Boussinesq-type system given in \eqref{Dee}.

\begin{thm}\label{MA3}
Let $L, T>0$  be fixed constants. Consider the initial data $N_0$ and $V_0$ $\in H_{0}^1$, and a given coefficient $c \in L^2\left[0,L \right]$. Then, there exists a unique solution $(N,V) \in \mathcal{L}_{T_0}^2 \times \mathcal{L}_{T_0}^2$ to the system $\left(\ref{Dee}\right)$ subject to conditions $(\ref{Maa})$ and $(\ref{BBMLin33})$. Moreover, this solution depends continuously on the initial data $(N_0,V_0)$ and satisfies the following estimate:
\begin{equation}\label{D1}
		\left\| (N,V)\right\|_{\mathcal{L}_{T_0}^2\times \mathcal{L}_{T_0}^2}\leq \left\| (N_0,V_0)\right\|_{H_{0}^1\times H_{0}^1} T_0^{\frac{1}{2}} \text{exp} \left[\frac{L^{\frac{1}{2}}}{\left( \frac{\beta}{6}\right)^{\frac{1}{2}} }D T_0 \right].
	\end{equation}
\end{thm}

\textbf{Proof:} By starting from the integral equations $\left(\ref{Ma11new}\right)$, $\left(\ref{Ma1new}\right)$, and taking into account the boundary conditions for $N$, $V$, we integrate over the interval $(0,t)$, to obtain
\begin{equation}\label{Ana11}
	N(\xi,t)= N_0(\xi)+\int_{0}^{t}\int_{0}^{L}K(\xi,s)\left(1+ \alpha c^{2}N(s,\tau) \right)  V(s,\tau)ds d\tau,
\end{equation}

\begin{equation}\label{Ana12}
\begin{split}
	V(\xi,t)= V_0(\xi) & +\int_{0}^{t}\int_{0}^{L}K(\xi,s)c(s)N(s,\tau)ds d\tau 
    \\ \\ & + \dfrac{1}{2}\int_{0}^{t}\int_{0}^{L}K(\xi,s)\alpha c^{2}(s)V^{2}(s,\tau)ds d\tau,
    \end{split}
\end{equation}
or equivalently,
\begin{multline}\label{Ana16}
	\left( N(\xi,t),V(\xi,t)\right) = \left( N_0(\xi) +\int_{0}^{t}\int_{0}^{L}K(\xi,s)\left(1+ \alpha c^{2}N(s,\tau) \right)  V(s,\tau)ds d\tau  \right.,\\ \left. V_0(\xi)+\int_{0}^{t}\int_{0}^{L}K(\xi,s)cN(s,\tau)ds d\tau 
     + \dfrac{1}{2}\int_{0}^{t}\int_{0}^{L}K(\xi,s)\alpha c^{2}V^{2}(s,\tau)ds d\tau  \right), 
\end{multline}
where $G(\xi,s)$ and $K(\xi,s)$ are the functions defined in \eqref{FuncionG}, \eqref{DerivadaG}. Let us consider the nonlinear operator
\begin{equation}
	\mathcal{A}: \mathcal{L}_{T}^2 \times \mathcal{L}_{T}^2 \rightarrow \mathcal{L}_{T}^2 \times \mathcal{L}_{T}^2,
\end{equation}
defined by
\begin{equation}\begin{array}{c l}\label{Ma2}
		\mathcal{A}\left(N(\xi,t), V(\xi,t) \right)=\\\\\left( N_0(\xi)+ \int_{0}^{t}\int_{0}^{L}K(\xi,s) V(s,\tau)ds d\tau  +\int_{0}^{t}\int_{0}^{L}K(\xi,s) \alpha c^{2}N(s,\tau)  V(s,\tau)ds d\tau \right.,\\
		\left. V_0(\xi)+\int_{0}^{t}\int_{0}^{L}K(\xi,s)N(s,\tau)ds d\tau + \dfrac{1}{2}\int_{0}^{t}\int_{0}^{L}K(\xi,s)\alpha c^{2}V^{2}(s,\tau)ds d\tau  \right).
	\end{array} 
\end{equation}
Note that from the equalities $\left(\ref{Ana11}\right)$,$\left(\ref{Ana12}\right)$, it can be established that
$$\mathcal{A}\left(N(0,t),V(0,t) \right) \equiv (0,0),  ~~ \mathcal{A}\left(N(L,t),V(L,t)\right) \equiv (0,0),$$ 
for all $t\in \left[0,T \right] $.

Therefore, establishing the existence of solutions to the integral equations $\left(\ref{Ana11}\right)$ and $\left(\ref{Ana12}\right)$, is equivalent to proving the existence of a fixed point for the operator $\mathcal{A}$ in the product space $\mathcal{L}_{T_0}^2 \times \mathcal{L}_{T_0}^2$, where $T_0>0$ is chosen appropriately. 

In order to apply the Banach fixed point theorem, we consider the closed ball of radius $R$
\begin{equation}
	\mathcal{B}_{R_T}=\left\lbrace (N,V) \in \mathcal{L}_{T}^2 \times \mathcal{L}_{T}^2:  \left\| (N,V)\right\|_{\mathcal{L}_{T}^2 \times \mathcal{L}_{T}^2} \leq R\right\rbrace,
\end{equation}
and let $(N,V) \in \mathcal{B}_{R_{T_{1}}}$, where $0< T_{1}\leq T$. Note that from the definition $\left(\ref{Ma2}\right)$ and using Lemma (\ref{MA51}) it follows that

\begin{multline*}
	\left\| \mathcal{A}\left(N(\xi,t),V(\xi,t)\right) \right\|_{H_{0}^1\times H_{0}^1}  \leq \left\|N_0 \right\| _{H_{0}^1}  +\left\| \int_{0}^{t}\int_{0}^{L}K(\cdot ,s) V(s,\tau)ds d\tau   \right\|_{H_{0}^1}  \\\\ +\left\|  \int_{0}^{t}\int_{0}^{L}K( \cdot ,s) \alpha c^{2}(s)N(s,\tau)  V(s,\tau)ds d\tau  \right\|_{H_{0}^1} 
	+ \left\|V_0 \right\| _{H_{0}^1} \\\\ +\left\|\int_{0}^{t}\int_{0}^{L}K(\cdot,s)c(s) N ds d\tau \right\|_{H_{0}^1}  + \left\|  \dfrac{1}{2}\int_{0}^{t}\int_{0}^{L}K(\cdot,s)\alpha c^{2}(s)V^{2}(s,\tau)ds d\tau   \right\|_{H_{0}^1}.
\end{multline*}
Consequently,
\begin{multline*}
	\left\| \mathcal{A}\left(N(\xi,t),V(\xi,t)\right) \right\|_{H_{0}^1\times H_{0}^1}  \leq  \left\|N_0 \right\| _{H_{0}^1}+ \left\|V_0 \right\| _{H_{0}^1} \\\\  + \int_{0}^{t}\left\| \int_{0}^{L}K(\cdot ,s) V(s,\tau)ds\right\|_{H_{0}^1} d\tau   +  \int_{0}^{t}\left\|  \int_{0}^{L}K( \cdot ,s) \alpha c^{2}(s)N(s,\tau)  V(s,\tau)ds \right\|_{H_{0}^1}d\tau 
	\\ \\ +\int_{0}^{t} \left\| \int_{0}^{L}K(\cdot,s)c(s) N ds \right\|_{H_{0}^1} d\tau   +  \dfrac{1}{2}\int_{0}^{t}\left\|  \int_{0}^{L}K(\cdot,s)\alpha c^{2}(s)V^{2}(s,\tau)ds  \right\|_{H_{0}^1}d\tau .
\end{multline*}
Therefore, using the Lemma \ref{MA51}, we get
\begin{multline*}
	\left\| \mathcal{A}\left(N(\xi,t),V(\xi,t)\right) \right\|_{H_{0}^1\times H_{0}^1}  \leq \left\|N_0 \right\|_{H_0^1}+\left\|V_0 \right\| _{H_{0}^1}+\int_{0}^{t}C_1\left\|V \right\|_{L^{2}}d\tau \\ +\alpha\int_{0}^{t}C_2\left\| c^{2} N  V \right\|_{L^2}d\tau
	+\int_{0}^{t}C_3\left\|cN \right\|_{L^2}d\tau+ \dfrac{\alpha}{2}\int_{0}^{t}C_4 \left\| c^{2}V^{2}\right\|_{L^2}d\tau.
\end{multline*}
Taking $C_1=C_3$, $C_2=C_4$ and using the continuous inclusion of $H_{0}^1\left(\left[0,L\right] \right) $ in $L^{\infty}\left( \left[0, L\right]\right) $, we have
\begin{multline*}
	\left\| \mathcal{A}\left(N(\xi,t),V(\xi,t)\right) \right\|_{H_{0}^1\times H_{0}^1}  \\ \leq \left\|N_0 \right\|_{H_0^1}+\left\|V_0 \right\| _{H_{0}^1}  +C_1\int_{0}^{t} L ^{\frac{1}{2}}\left\| V\right\| _{L^{\infty}} + \left\|c \right\|_{L^2}\left\|N\right\| _{L^{\infty}} d \tau \\  + \alpha \: C_2 \left\|c^{2} \right\|_{L^2}  \int_{0}^{t}\left\| N\right\| _{L^{\infty}} \left\| V\right\| _{L^{\infty}}+  \dfrac{1}{2} \left\| V\right\|^{2} _{L^{\infty}} d\tau.
\end{multline*}
Therefore, we obtain that
\begin{multline*}
	\left\| \mathcal{A}\left(N(\xi,t),V(\xi,t)\right) \right\|_{H_{0}^1\times H_{0}^1}  \\ \leq \left\|N_0 \right\|_{H_0^1}+\left\|V_0 \right\| _{H_{0}^1}+\frac{C_1L^{\frac{1}{2}}}{\left(\frac{\beta}{6}\right)^{\frac{1}{2}} }\int_{0}^{T_1}  L^{\frac{1}{2}} \left\| V\right\| _{H_{0}^{1}}+ \left\|c \right\|_{L^2}\left\|N\right\| _{H_{0}^{1}} d\tau  \\  + \dfrac{\alpha \: C_2 \left\|c^{2} \right\|_{L^2} L}{\left( \dfrac{\beta}{6}\right) } \int_{0}^{T_1} \left\| N\right\| _{H_{0}^{1}}  \left\| V\right\| _{H_{0}^{1}}+ \dfrac{1}{2} \left\| V\right\|^{2} _{H_{0}^{1}}   d\tau.
\end{multline*}
The Cauchy-Schwarz inequality implies that
\begin{multline*}
	\left\| \mathcal{A}\left(N(\xi,t),V(\xi,t)\right) \right\|_{H_{0}^1\times H_{0}^1}  \\  \leq \left\|N_0 \right\|_{H_0^1}+\left\|V_0 \right\| _{H_{0}^1}+ \frac{C_1T_{1}^{\frac{1}{2}}L^{\frac{1}{2}}}{\left(\frac{\beta}{6}\right)^{\frac{1}{2}} }\left(L^{\frac{1}{2}} \left\| V\right\|_{\mathcal{L}_{T_1}^{2}}+ \left\|c \right\|_{L^2}\left\|N\right\| _{\mathcal{L}_{T_1}^{2}}\right)  \\  +\dfrac{\alpha \: C_2 \left\|c^{2} \right\|_{L^2} L T_{1}}{\left( \dfrac{\beta}{6}\right) }\left( \left\| N\right\|_{\mathcal{L}_{T_1}^{2}} \left\| V\right\|_{\mathcal{L}_{T_1}^{2}}+ \dfrac{1}{2} \left\| V\right\|^{2}_{\mathcal{L}_{T_1}^{2}}\right).
\end{multline*}
Then, since $(N,V)\in B_{R_{T_{1}}}$ , it follows that

\begin{multline*}
	\left\| \mathcal{A}\left(N(\xi,t),V(\xi,t)\right) \right\|_{H_{0}^1\times H_{0}^1}  \\ \leq \left\|N_0 \right\|_{H_0^1}+\left\|V_0 \right\| _{H_{0}^1}+ \frac{C_1T_{1}^{\frac{1}{2}}L^{\frac{1}{2}}}{\left(\frac{\beta}{6}\right)^{\frac{1}{2}} }\left(L^{\frac{1}{2}} R+ \left\|c \right\|_{L^2}R \right)   \\  +\dfrac{\alpha \: C_2 \left\|c^{2} \right\|_{L^2} L T_{1}}{\left( \dfrac{\beta}{6}\right) }\left(  R^{2}+ \dfrac{1}{2} R^{2}\right).
\end{multline*}
Squaring and integrating over $[0,T_1]$, we obtain

\begin{multline*}
	\left\| \mathcal{A}\left(N(\xi,t),V(\xi,t)\right) \right\|_{\mathcal{L}_{T_1}^{2}\times \mathcal{L}_{T_1}^{2}}  \\\leq T_{1}^{\frac{1}{2}} \left[  \left\|N_0 \right\|_{H_0^1}+\left\|V_0 \right\| _{H_{0}^1}+ \frac{C_1T_{1}^{\frac{1}{2}}L^{\frac{1}{2}}R}{\left(\frac{\beta}{6}\right)^{\frac{1}{2}} }\left(L^{\frac{1}{2}} + \left\|c \right\|_{L^2} \right)\right]   \\ +T_{1}^{\frac{1}{2}}\left[ \dfrac{3\alpha \: C_2 \left\|c^{2} \right\|_{L^2} L T_{1}^{\frac{1}{2}}R^{2}}{2\left( \dfrac{\beta}{6}\right) }  \right] .
\end{multline*}
Then, if we take $T_1 > 0$ fixed and small enough such that


\begin{multline}\label{Ana13}
	T_{1}^{\frac{1}{2}} \left[ \left\|N_0 \right\|_{H_{0}^1}+\left\|V_0 \right\| _{H_{0}^1} \right. \\  + \left.\frac{L^{\frac{1}{2}}}{\left(\frac{\beta}{6}\right)^{\frac{1}{2}} }\left[ C_1 \left( L^{\frac{1}{2}}+ \left\|c \right\|_{L^2}\right)+\dfrac{3\alpha \: C_2 \left\|c^{2} \right\|_{L^2} L ^{\frac{1}{2}}T_{1}^{\frac{1}{2}}R}{2\left( \dfrac{\beta}{6}\right)^{\frac{1}{2}} }  \right]T_{1}^{\frac{1}{2}}R  \right] 
    \leq R,
\end{multline}
we get
\begin{equation*}
	\left\| \mathcal{A}\left(N(\xi,t),V(\xi,t)\right) \right\|_{\mathcal{L}_{T_1}^{2}\times \mathcal{L}_{T_1}^{2}} \leq R.
\end{equation*}
It follows that the operator $\mathcal{A}$ maps the ball into itself, that is,
$$\mathcal{A}\left( B_{R_{T_{1}}}\right)\subset B_{R_{T_{1}}},$$ 
for sufficiently small $T_1 > 0$.
Next, to establish conditions under which the operator $\mathcal{A}$ 
is a contraction on $B_{R_{T_{2}}}$, with
$0 < T_{2} \leq T_1$ chosen appropriately, consider two arbitrary elements $\left( N,V\right) $ and $\left( N_{*},V_{*}\right)$ in $B_{R_{T_{2}}}$. Observe that
\begin{multline*}
	( \mathcal{A}(N,V)-\mathcal{A}(N_*,V_*))(\xi,t)= \left(\int_0^t\int_{0}^{L} K(\xi,s)(V-V_* )dsd\tau \right.  \\
	\left. +\alpha \int_0^t\int_{0}^{L} K(\xi,s)c^{2}(s)(NV-N_*V_*) dsd\tau, 
	\int_0^t\int_{0}^{L} \! K(\xi,s) c(s)(N-N_*) dsd\tau   \right.  \\
	\left.  + \dfrac{\alpha}{2}	\int_0^t\int_{0}^{L} \! K(\xi,s) c^{2}(s)(V^{2}-V^2_*) dsd\tau \right) .
\end{multline*}
By applying an analogous procedure as before, we establish that
\begin{multline*}
	\left\|\mathcal{A}(N,V)-\mathcal{A}(N_*,V_*)\right\|_{\mathcal{L}_{T_{2}}^{2}\times \mathcal{L}_{T_{2}}^{2}}\\ \leq \frac{L^{\frac{1}{2}}T_{2}}{\left(\frac{\beta}{6} \right)^{\frac{1}{2}} }\left[C_1\left(L^{\frac{1}{2}}+\left\|c \right\|_{L^{2}} \right)+\alpha \: C_2\left\| c^{2} \right\|_{L^{2}}\left( 1+2R\right) \right] \left\|(N,V)-(N_*,V_*)\right\|_{\mathcal{L}_{T_{2}}^{2}\times \mathcal{L}_{T_{2}}^{2}}.
\end{multline*}
Therefore, if we take $T_2 > 0$ fixed and sufficiently small, such that
\begin{equation}\label{Ana15}
	\frac{L^{\frac{1}{2}}T_{2}}{\left(\frac{\beta}{6} \right)^{\frac{1}{2}} }\left[C_1\left(L^{\frac{1}{2}}+\left\|c \right\|_{L^{2}} \right)+\alpha \: C_2\left\| c^{2} \right\|_{L^{2}}\left( 1+2R\right) \right]	< 1,
\end{equation}

we have that the operator $\mathcal{A}$ is a contraction in $B_{R_{T_{2}}}$. It is important to note that the constants that
accompany $T_1$ and $T_2$ in $\left(\ref{Ana13}\right)$, and $\left(\ref{Ana15}\right)$, are independent of $T_1$ and $T_2$, respectively.
Finally, consider $T_{0} = \min\left\lbrace T_{1} , T_{2} \right\rbrace$. As a consequence of the Banach Fixed Point Theorem applied to $B_{R_{T_{0}}}$ , the existence of a unique $(N,V) \in B_{R_{T_{0}}}$ satisfying the integral equation $\left(\ref{Ana16}\right)$ is established.\\

Furthermore, to establish the uniqueness of the solutions to the \eqref{Dee}-\eqref{BBMLin33} problem in $\mathcal{L}_{T_{0}}^{2}\times \mathcal{L}_{T_{0}}^{2} $ , consider
$(N,V)$ and $(\tilde{N},\tilde{V})$ elements in $\mathcal{L}_{T_{0}}^{2}\times \mathcal{L}_{T_{0}}^{2} $ solutions to the integral equation $\left(\ref{Ana16}\right)$ with initial data $(N_{0},V_{0})$ and $(\tilde{N}_{0},\tilde{V}_{0})$ , 
respectively. Note that

\begin{multline*}
	((N,V)-(\tilde{N},\tilde{V}))(\xi,t)=\\  \left( \left(N_{0}-\tilde{N}_{0}\right) + \int_0^t\int_{0}^{L} K \left( V-\tilde{V}\right)  dsd\tau \right.
	+ \int_0^t\int_{0}^{L} K \alpha c ^{2}(NV-\tilde{N}\tilde{V}) dsd\tau,\\ \left(V_{0}-\tilde{V}_{0}\right)+\left. \int_0^t\int_{0}^{L} \! K c (N-\tilde{N}) dsd\tau + \frac{1}{2} \int_0^t\int_{0}^{L} \! K \alpha c^{2} (V^{2}-\tilde{V}^{2}) dsd\tau   \right).
\end{multline*}
A direct calculation, similar to the previous ones, allows us to establish that

\begin{multline*}
	\left\|(N,V)-(\tilde{N},\tilde{V})\right\|_{H_{0}^{1}\times H_{0}^{1} } \\ \leq  \left\|(N_0,V_0)-(\tilde{N}_0,\tilde{V}_0)\right\|_{H_{0}^{1}\times H_{0}^{1} }+ \frac{L^{\frac{1}{2}}}{\left(\frac{\beta}{6} \right)^{\frac{1}{2}} }\int_{0}^{t}\left[C_1\left(L^{\frac{1}{2}}+\left\|c \right\|_{L^{2}} \right)\right. \\ \left. +\alpha \: C_2\left\| c^{2} \right\|_{L^{2}}\left( 2\left\| V \right\|_{L_{T_{2}}}+\left\| \tilde{N} \right\|_{L_{T_{2}}}+ \left\| \tilde{V} \right\|_{L_{T_{2}}}\right) \right] \left\|(N,V)-(\tilde{N},\tilde{V})\right\|_{H_{0}^{1}\times H_{0}^{1} }d\tau.
\end{multline*}
In view of the above, Gronwall's lemma \cite{Evans} implies that
\begin{multline*}
	\left\|(N,V)-(\tilde{N},\tilde{V})\right\|_{H_{0}^{1}\times H_{0}^{1} } \\ \leq  \left\|(N_0,V_0)-(\tilde{N}_0,\tilde{V}_0)\right\|_{H_{0}^{1}\times H_{0}^{1} } \exp \left[\frac{L^{\frac{1}{2}}}{\left(\frac{\beta}{6} \right)^{\frac{1}{2}} }\int_{0}^{t}\left(C_1\left(L^{\frac{1}{2}}+\left\|c \right\|_{L^{2}} \right)\right.\right.\\ + \left. \left.\alpha \: C_2\left\| c^{2} \right\|_{L^{2}}\left( 2\left\| V \right\|_{L_{T_{2}}}+\left\| \tilde{N} \right\|_{L_{T_{2}}}+ \left\| \tilde{V} \right\|_{L_{T_{2}}}\right) \right)  d\tau\right] \\
	\leq \left\|(N_0,V_0)-(\tilde{N}_0,\tilde{V}_0)\right\|_{H_{0}^{1}\times H_{0}^{1} } \exp\left[\frac{L^{\frac{1}{2}}}{\left(\frac{\beta}{6} \right)^{\frac{1}{2}} }\int_{0}^{T_0}\left(C_1\left(L^{\frac{1}{2}}+\left\|c \right\|_{L^{2}} \right)\right.\right.\\ + \left.\left. \alpha \: C_2\left\| c^{2} \right\|_{L^{2}}\left(2\left\| V \right\|_{L_{T_{2}}}+\left\| \tilde{N} \right\|_{L_{T_{2}}}+ \left\| \tilde{V} \right\|_{L_{T_{2}}}\right)  \right)  d\tau\right].
\end{multline*}
Squaring and integrating in $\left(0,T_0 \right)$, we arrive at

\begin{multline*}
	\left\|(N,V)-(\tilde{N},\tilde{V})\right\|_{\mathcal{L}_{T_0}^{2}\times \mathcal{L}_{T_0}^{2} } \leq \\
	\left\|(N_0,V_0)-(\tilde{N}_0,\tilde{V}_0)\right\|_{H_{0}^{1}\times H_{0}^{1} }T_{0}^{\frac{1}{2}}\exp \left[  \frac{L^{\frac{1}{2}}}{\left(\frac{\beta}{6} \right)^{\frac{1}{2}} }\left(C_1\left(L^{\frac{1}{2}}+\left\|c \right\|_{L^{2}} \right)\right.\right. \\ \left.\left. +\alpha \: C_2\left\| c^{2} \right\|_{L^{2}}\left( 2\left\| V \right\|_{L_{T_{2}}}+\left\| \tilde{N} \right\|_{L_{T_{2}}}+ \left\| \tilde{V} \right\|_{L_{T_{2}}}\right) \right)T_0 \right].
\end{multline*}
Therefore, it can be established that

\begin{multline*}
	\left\|(N,V)-(\tilde{N},\tilde{V})\right\|_{\mathcal{L}_{T_0}^{2}\times \mathcal{L}_{T_0}^{2} }  \leq 
	\left\|(N_0,V_0)-(\tilde{N}_0,\tilde{V}_0)\right\|_{H_{0}^{1}\times H_{0}^{1} }T_{0}^{\frac{1}{2}}\exp \left[  \frac{L^{\frac{1}{2}}}{\left(\frac{\beta}{6} \right)^{\frac{1}{2}} }D T_0 \right].
\end{multline*}
The previous estimate establishes the continuous dependence on the initial data and also the uniqueness of the solutions to the integral equation $\left(\ref{Ana16}\right)$ in $\mathcal{L}^{2}_{T_0} \times \mathcal{L}^{2}_{T_0}$.
It has been established that there exists a unique solution to the problem defined in $[0,T_0]$ for $T_0$
sufficiently small.
$\hfill\square$

\begin{thm}\label{MA3}
	Let $L>0$, and let $\alpha, \beta >0$ be fixed constants. Suppose the initial data $N_0$ and $V_0$ $\in L^2\left[0,L \right]$ and let the coefficient function $c \in L^2\left[0,L \right]$ be given. Then any solution $(N,V) \in \mathcal{L}_{T_0}^2 \times \mathcal{L}_{T_0}^2$ to the system $\left(\ref{Dee}\right)$, subject to the initial and boundary conditions $(\ref{Maa})$ and $(\ref{BBMLin33})$, satisfies the following conservation law:
	\begin{equation}\label{D1}
		\left\| V \right\|^{2}_{L^2 }+ \alpha\left\| c |N|^{\frac{1}{2}}V \right\|^{2}_{L^2 } +\left\| \left| c\right| ^{\frac{1}{2}}N  \right\|^{2}_{L^2 }  = \left\| V_0\right\|^{2}_{L^2 } +  \alpha\left\| c |N_0|^{\frac{1}{2}}V_0 \right\|^{2}_{\mathcal{L}^2 }+ \left\| \left| c\right| ^{\frac{1}{2}}N_0 \right\|^{2}_{L^2}.
	\end{equation}
In other words, the energy functional $E(t)$ is conserved in time, i.e. 
$E(t)=E(0)$ for all $t$.
\end{thm}

\textbf{Proof:} 
Let 
\begin{multline*}
	E(t) = \frac{1}{2} \left\| V \right\|^{2}_{L^2 }+ \frac{\alpha}{2} \left\| c (N)^{\frac{1}{2}}V \right\|^{2}_{L^2} + \frac12 \left\| \left| c\right|^{\frac{1}{2}}N \right\|^{2}_{L^2 }  = \frac{1}{2}\int_{0}^{L}\left[ \left(1+\alpha c^{2}N\right) V^{2}+cN^{2}\right].
    \end{multline*}
Differentiating $E(t)$, we have
\begin{equation}\label{D2}
	\dfrac{dE(t)}{dt}=\int_{0}^{L}\left[  \alpha c^{2}N_t V^{2}+(1+\alpha c^{2}N)VV_{t}+c NN_t\right] d\xi,
\end{equation}
and by adding and subtracting appropriate terms, we obtain
\begin{multline*}
	\dfrac{dE(t)}{dt}=\int_{0}^{L}\left[(1+\alpha c^{2}N)V-\dfrac{\beta}{6}\partial_\xi N_t\right]V_t\: d\xi \\+ \int_{0}^{L}\left[cN + \frac12 \alpha c^{2}V^{2}-\dfrac{\beta}{6}\partial_\xi V_t\right]N_t\: d\xi   + \dfrac{\beta}{6}\int_{0}^{L} \partial_\xi \left(N_t V_t \right)d\xi.
\end{multline*}
Note that
\begin{equation*}
	\int_{0}^{L} \partial_\xi \left(N_t V_t \right)d\xi = N_t\left(L,t \right)V_t\left(L,t \right)- N_t\left(0,t \right)V_t\left(0,t \right)=0.
\end{equation*}
due to $N \left(L,t \right)= N \left(0,t \right)=0$, for all $t$. Therefore, we deduce that
\begin{equation}
	\dfrac{dE(t)}{dt}=\int_{0}^{L}\partial \xi\left[\left( \left( 1+\alpha c^{2}N\right)V -\dfrac{\beta}{6} N_{\xi t} \right)\left(cN+\dfrac{\alpha}{2}c^{2}V^{2}-\dfrac{\beta}{6} V_{\xi t} \right)  \right] d\xi=0, 
\end{equation}
In other words, $E(t)=E(0)$, that is,
\begin{equation}
	\left\| V \right\|^{2}_{ L^2 }+ \alpha\left\| c (N)^{\frac{1}{2}}V \right\|^{2}_{L^2 } +\left\| \left| c\right| ^{\frac{1}{2}}N  \right\|^{2}_{L^2 }  = \left\| V_0\right\|^{2}_{L^2 } +  \alpha\left\| c (N_0)^{\frac{1}{2}}V_0 \right\|^{2}_{L^2}+ \left\| \left| c\right| ^{\frac{1}{2}}N_0  \right\|^{2}_{L^2}, 
\end{equation}
for all $t$. $\hfill\square$

\section{Numerical experiments}

In this section, we perform some numerical simulations with the nonlinear Boussinesq-type system
\begin{equation}\label{MA29}
\begin{array}{c l}
N_t+\left[  \left( 1+ \alpha c^{2}N\right) V\right] _{\xi}-\frac{\beta}{6}N_{\xi\xi t}&=0 \:\:\: ,\:\:\left(\xi,t\right)\in \left[0,L\right]\times\left[0,T\right].\\
V_t+(cN)_{\xi}+\frac{1}{2} \alpha \left( c^{2}V^{2}\right)_{\xi} -\frac{\beta}{6}V_{\xi\xi t}&=0,
\end{array}
\end{equation}
subject to the initial conditions
\begin{equation}\label{datosiniciales}
\begin{array}{c l}
N(\xi,0)&=N_0(\xi),\\
V(\xi,0)&=V_0(\xi),\:\:\: \xi \in \left[0,L\right],
\end{array}
\end{equation}
and Dirichlet-type boundary conditions
\begin{equation}\label{condfrontera}
\begin{array}{c l}
N(0,t)=N(L,t)=0,\:\: t\in \left[0,T\right]\\
V(0,t)=V(L,t)=0,\:\: t\in \left[0,T\right].
\end{array}
\end{equation}

To this end, we use the finite element-space numerical scheme 
proposed in \cite{PizoMunoz}, which we recall that can be written as



\begin{equation}\label{sistemaNV1}
\partial_t\left( \left\langle N,v_1\right\rangle _{L^{2}}+\frac{\beta}{6}\left\langle N_{\xi },v_{1\xi}\right\rangle _{L^{2}}\right) =\left\langle\left[  \left( 1 + \alpha c^{2}N\right)V\right] ,v_{1\xi}\right\rangle _{L^{2}}
\end{equation}
\begin{equation}\label{sistemaNV2}
\partial_t\left( \left\langle V,v_2\right\rangle _{L^{2}}+\frac{\beta}{6}\left\langle V_{\xi },v_{2\xi}\right\rangle _{L^{2}}\right) =\left\langle (cN),v_{2\xi}\right\rangle_{L^{2}} +\frac{1}{2}\alpha \left\langle c^{2}V^{2} ,v_{2\xi}\right\rangle _{L^{2}},
\end{equation}
where $v_1 , v_2 \in H_0^1(0,L)$ are test functions.
System \eqref{sistemaNV1}-\eqref{sistemaNV2} is discretized in the temporal variable by the following numerical scheme:
\begin{multline}\label{MA33}
 \left\langle \frac{N^{n+1}-N^{n}}{\Delta t},v_1\right\rangle _{L^{2}}+\frac{\beta}{6}\left\langle \frac{N_{\xi }^{n+1}-N_{\xi}^{n}}{\Delta t},v_{1\xi}\right\rangle _{L^{2}} \\\\ =\theta \left\langle  \left( 1+\alpha c^{2}N^{n+1}\right)V^{n+1} ,v_{1\xi}\right\rangle _{L^{2}} +(1-\theta)\left\langle \left( 1+ \alpha c^{2}N^{n}\right)V^{n},v_{1\xi}\right\rangle 
\end{multline}
\begin{multline}\label{MA34}
 \left\langle \frac{V^{n+1}-V^{n}}{\Delta t},v_2\right\rangle _{L^{2}}+\frac{\beta}{6}\left\langle \frac{V^{n+1}_{\xi }-V^{n}_{\xi}}{\Delta t},v_{2\xi}\right\rangle _{L^{2}}\\\\ =\theta \left\langle c N^{n+1},v_{2\xi}\right\rangle_{L^{2}}+\theta \frac{1}{2} \alpha \left\langle c^{2}V^{n+1}V^{n+1} ,v_{2\xi}\right\rangle _{L^{2}}+(1-\theta)\left\langle cN^{n},v_{2\xi}\right\rangle_{L^{2}}\\ +(1-\theta)\frac{1}{2} \alpha \left\langle c^{2}V^{n}V^{n} ,v_{2\xi}\right\rangle _{L^{2}},
\end{multline}
where $(N^{0},V^{0})$ is the initial condition $(N_0(\xi),V_0(\xi))$, and $(N^{n},V^{n})$ denotes the numerical approximation of the solution $(N,V)$ at time $t=n\Delta t$. The parameter $\theta$ is a real constant in the interval $(0, 1)$, which is set to $1/2$ in the experiments presented in this paper. Further details on the numerical solver can be found in \cite{PizoMunoz}. The numerical implementation was carried out  using the FEniCS library \cite{Fenics} in Python, which provides a flexible and efficient framework for implementing the finite element method to approximate solutions of initial-boundary value problems associated with partial differential equations.

As a first step, we compute numerical solutions of the initial-boundary problem associated with the variable-coefficient Boussinesq system \eqref{MA29}, subject to the initial conditions \eqref{datosiniciales} and the boundary conditions \eqref{condfrontera}.
The simulations are carried out for selected forms of the coefficient
function $c(\xi)$ and fixed values of the model parameters. The initial data are chosen as localized Gaussian profiles:
\[
N(\xi, 0) = N_0(\xi) = e^{-(\xi-18)^2}, ~~~ V(\xi,0) = V_0(\xi) =
e^{-(\xi-18)^2 },
\]
which correspond to a single smooth wave centered at $\xi = 18$. These profiles are symmetric and rapidly decaying, and are commonly used in dispersive wave simulations due to their smoothness and localized nature.

\begin{figure}[ht]
  \centering
	   \begin{subfigure}{0.47\linewidth}
		\includegraphics[width=1\linewidth, height=1.2\linewidth]{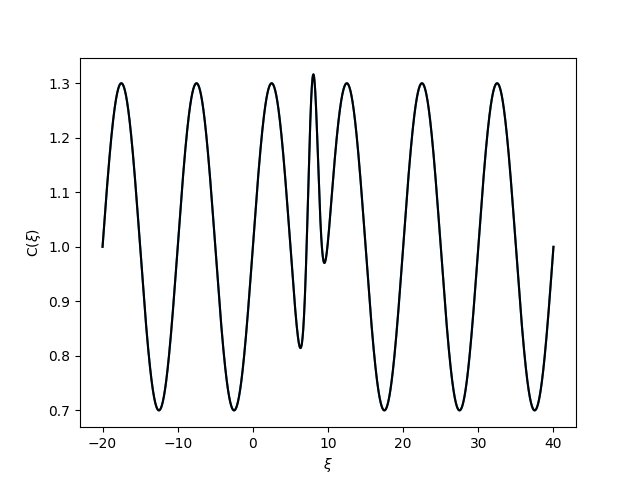}
		\caption{Coefficient $c(\xi)$ given in \eqref{coefficient_gauss}.}
	   \end{subfigure}
	   \begin{subfigure}{0.47\linewidth}
	      	\includegraphics[width=1\linewidth, height=1.2\linewidth]{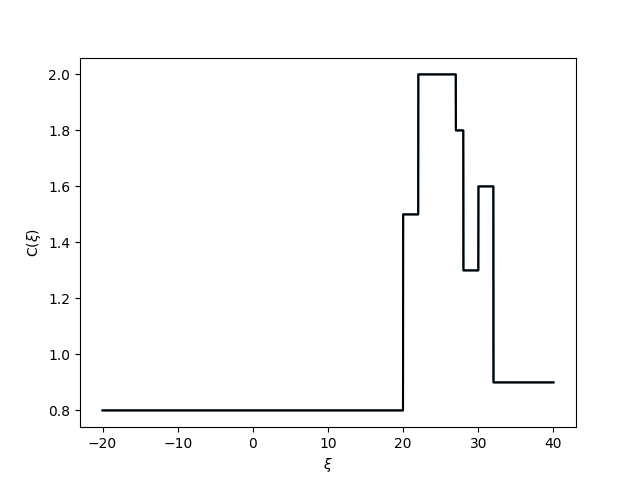}
		\caption{Coefficient $c(\xi)$ given in \eqref{coefficient_step}. }
	    \end{subfigure}
	\caption{Coefficients $c(\xi)$ used in the simulations.  }
	\label{coefficients}
\end{figure}

For the model parameters, we fix $\alpha=\beta = 0.1$, and we consider two distinct forms for the variable coefficient $c(\xi)$,
which represents the effect of the depth variation in the medium. The first choice is a smooth oscillatory-Gaussian perturbation (see Figure \ref{coefficients}(a) ) given by

\begin{equation}\label{coefficient_gauss}
c(\xi) = 1 + 0.3 \sin( (\pi/5) \xi ) + 0.6 e^{-(\xi -8)^2 },
\end{equation}
which introduces both periodic variation and a localized peak in the wave speed. The second choice is a discontinuous, piecewise-constant function simulating abrupt depth changes:
\begin{equation}\label{coefficient_step}
c(\xi) = \begin{cases}
            0.8, ~~ \xi < 20 \\
            1.5, ~~ 20< \xi \leq 22 \\
            2, ~~ 22 < \xi \leq 27 \\
            1.8, ~~ 27 < \xi \leq 28 \\
            1.3, ~~ 28 < \xi \leq 30 \\
            1.6, ~~ 30 < \xi \leq 32 \\
            0.9, ~~ 32 < \xi,
         \end{cases}
\end{equation}
designed to emulate wave propagation through a layered medium with variable bathymetry. 

\begin{figure}[ht]
  \centering
	   \begin{subfigure}{0.47\linewidth}
		\includegraphics[width=1\linewidth, height=1.2\linewidth]{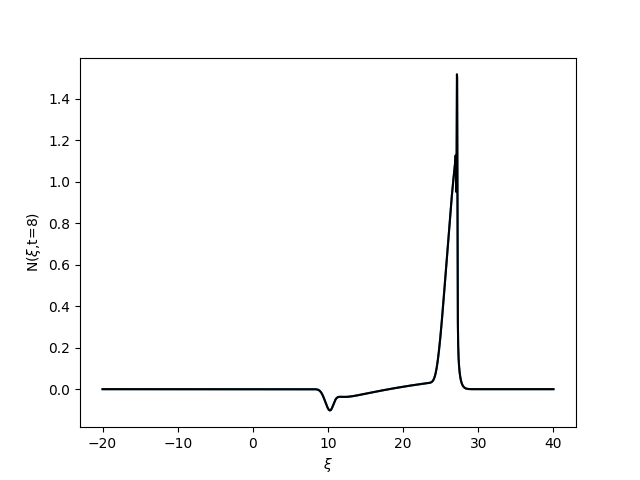}
		\label{N1}
	   \end{subfigure}
	   \begin{subfigure}{0.47\linewidth}
	      	\includegraphics[width=1\linewidth, height=1.2\linewidth]{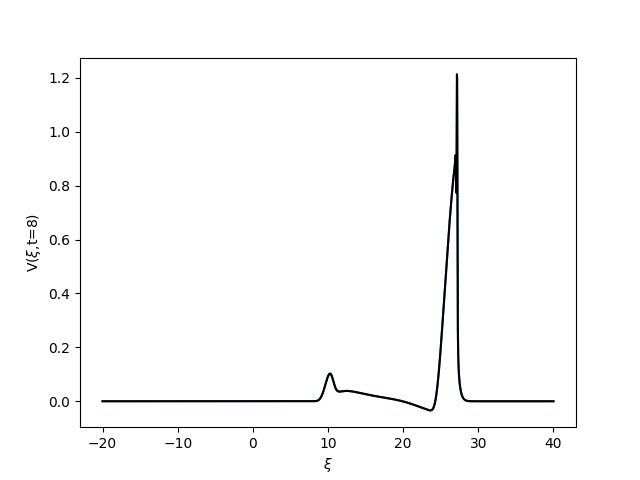}
		\label{V1}
	    \end{subfigure}
	\caption{Wave elevation $N(\xi,t)$ and fluid velocity $V(\xi,t)$ at $t = 8$, for the coefficient $c(\xi)$ shown in Figure \ref{coefficients}(a).  }
	\label{experiment_1}
\end{figure}

\begin{figure}[ht]
  \centering
	   \begin{subfigure}{0.47\linewidth}
		\includegraphics[width=1\linewidth, height=1.2\linewidth]{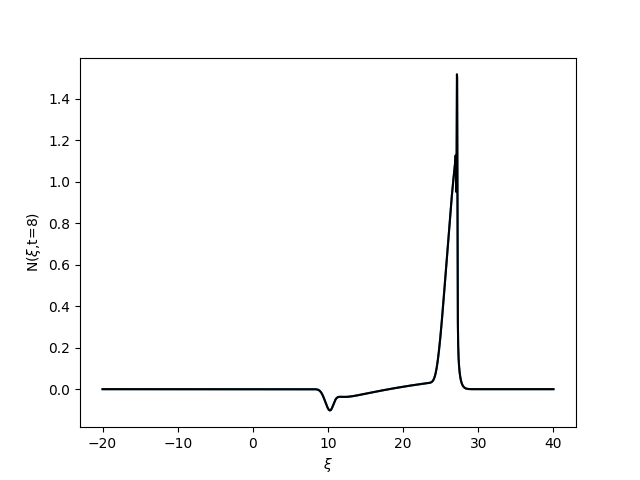}
		\label{N2}
	   \end{subfigure}
	   \begin{subfigure}{0.47\linewidth}
	      	\includegraphics[width=1\linewidth, height=1.2\linewidth]{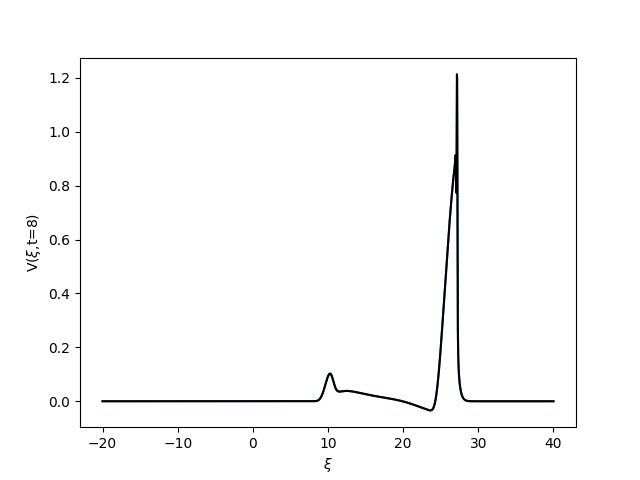}
		\label{V2}
	    \end{subfigure}
	\caption{Wave elevation $N(\xi,t)$ and fluid velocity $V(\xi,t)$ at $t = 8$ for the coefficient $c(\xi)$ shown in Figure \ref{coefficients}(b).  }
	\label{experiment_2}
\end{figure}

The numerical domain is the interval $[-20,40]$, discretized using
 3000 spatial grid points, resulting in a spatial resolution of $\Delta \xi \approx 0.02$. The time step is set to $\Delta t = 8/3000 \approx 2.7 e-3$, ensuring stability and accuracy of the time integration scheme. Figures \ref{experiment_1} and \ref{experiment_2} display the evolution of the numerical solution 
$(N(\xi,t), V(\xi,t)$ under the two different choices of the coefficient function $c(\xi)$. The results demonstrate how the variable depth influences the wave propagation, including changes in wave speed, amplitude modulation, and potential reflections or dispersive effects induced by the inhomogeneities in the medium.

These experiments serve as a foundation for more complex simulations and for the inverse problem discussed in subsequent sections, where the goal is to reconstruct initial conditions based on final-time observations.

Next, we illustrate the performance of the numerical solver \eqref{MA33}-\eqref{MA34} in solving an inverse problem associated with the variable-coefficient Boussinesq system. Given a final time $T$, model parameters $\alpha, \beta$, and a spatially-dependent coefficient $c(\xi)$, the goal is to reconstruct the initial conditions $N_0, V_0$, such that the solution $(N,V)$ to the initial-boundary value problem associated with system \eqref{MA29} satisfies the prescribed final-time observations
\[
N(\xi, T ) = N_T(\xi), ~~  V(\xi,T) = V_T(\xi), ~~~ \xi \in [0,L].
\]

The functions $N_T(\xi), V_T(\xi)$ are assumed to be known and are generated synthetically by solving the forward problem using the numerical scheme \eqref{MA33}–\eqref{MA34}, with prescribed initial conditions. This setup ensures consistency in the data and allows us to validate the reconstruction strategy.

To recover the unknown initial profiles $N_0$ and $V_0$, we formulate the inverse problem as an optimization task. Specifically, we define the cost functional 
\begin{equation}\label{Functional_assimilation}
J(N, V ) := \frac12 \int_0^L \Big( | N(\xi, T; N_0, V_0) - N_T(\xi) |^2 + |V(\xi,T; N_0, V_0) - V_T(\xi)|^2 \big) d\xi,
\end{equation}
that quantifies the discrepancy between the computed solution at time $T$ and the given observations. Here
$N(\xi, T; N_0, V_0), V(\xi, T; N_0, V_0)$ denote the numerical solution at time $T$, computed from the initial data $N_0, V_0$ using the forward solver.

The optimization problem consists in finding the initial pair
 $(N_0, V_0)$ that minimizes the functional $J$, thereby producing a solution trajectory that best matches the given observations at the final time. To formulate and solve this problem, we employ the FEniCS library \cite{Fenics} in combination with the Dolfin-Adjoint library \cite{dolfin_adjoint1, dolfin_adjoint2}, which provides automated adjoint capabilities and a convenient framework for defining and differentiating the objective functional $J$. The optimization is performed using the L-BFGS-B algorithm, a quasi-Newton method with limited memory and bound constraints, implemented via the optimization module in the SciPy library \cite{Nocedal}. This method is particularly well-suited for large-scale problems involving smooth functionals and offers efficient convergence even in high-dimensional settings.

\begin{figure}[ht]
  \centering
	   \begin{subfigure}{0.45\linewidth}
		\includegraphics[width=\linewidth, height=\linewidth]{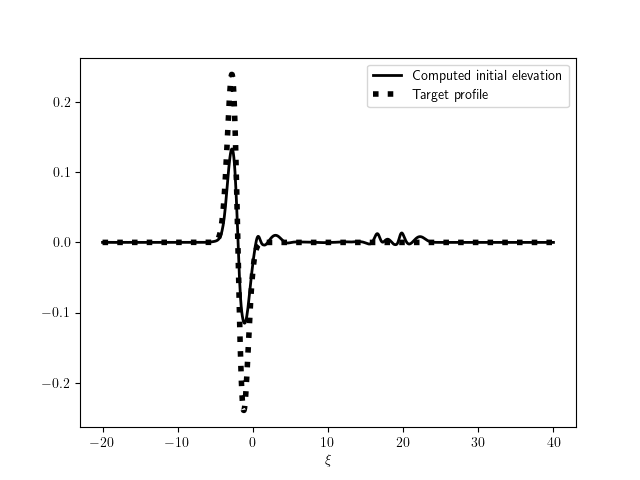}
		\caption{Iteration 1. }
	   \end{subfigure}
	   \begin{subfigure}{0.45\linewidth}
	      	\includegraphics[width=\linewidth, height=\linewidth]{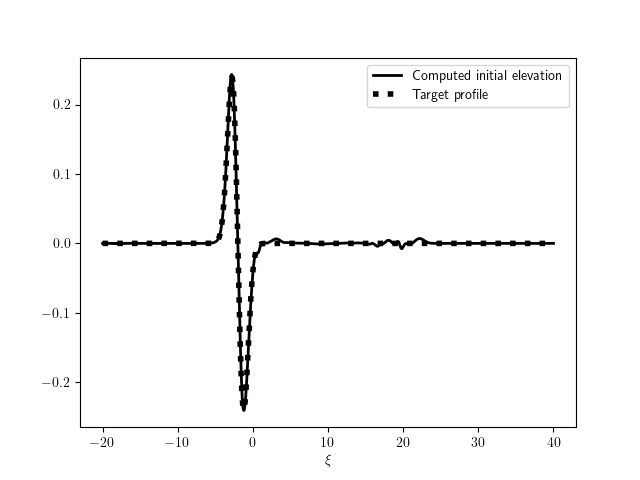}
		\caption{Iteration 2. }
	    \end{subfigure}
	\vfill
	     \begin{subfigure}{0.45\linewidth}
		 \includegraphics[width=\linewidth, height=\linewidth]{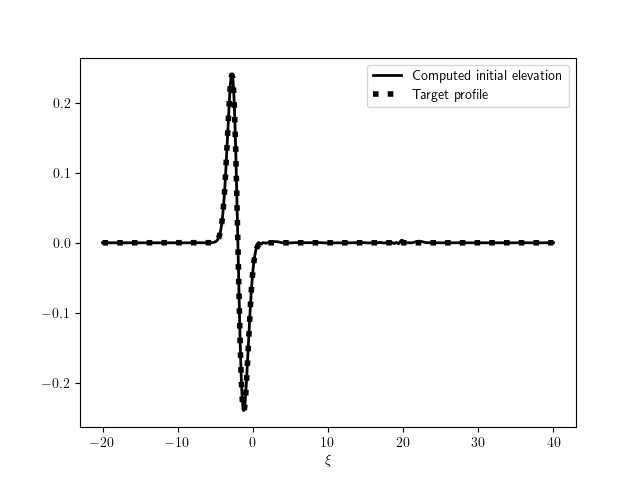}
		 \caption{Iteration 3. }
		 \label{}
	      \end{subfigure} 
	       \begin{subfigure}{0.45\linewidth}
		  \includegraphics[width=\linewidth, height=\linewidth ]{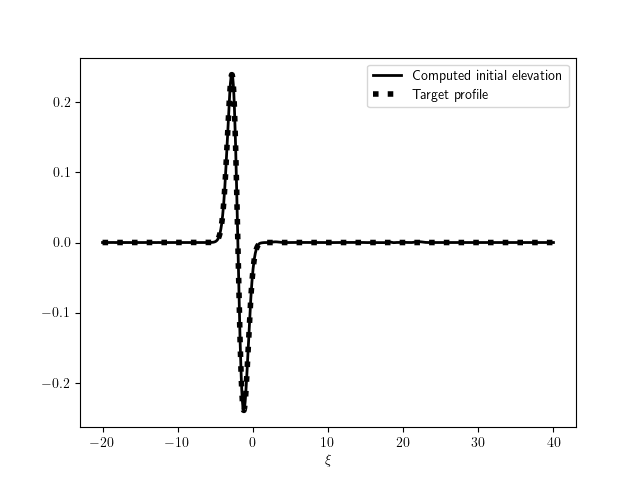}
		  \caption{Iteration 4. }
	       \end{subfigure}   
	\caption{Approximation of the initial wave elevation $N_0(\xi)$ after the first four iterations of the L-BFGS-B algorithm used to minimize the functional \eqref{Functional_assimilation}. }
	\label{Iterations_N}
\end{figure}

\begin{figure}[ht]
  \centering
	   \begin{subfigure}{0.45\linewidth}
		\includegraphics[width=\linewidth, height=\linewidth]{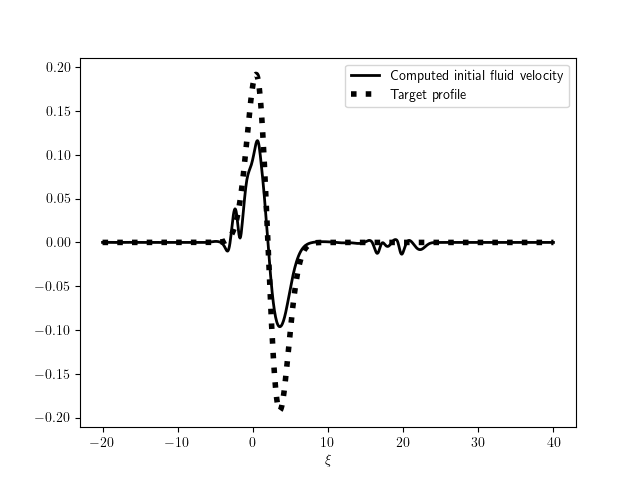}
		\caption{Iteration 1. }
	   \end{subfigure}
	   \begin{subfigure}{0.45\linewidth}
	      	\includegraphics[width=\linewidth, height=\linewidth]{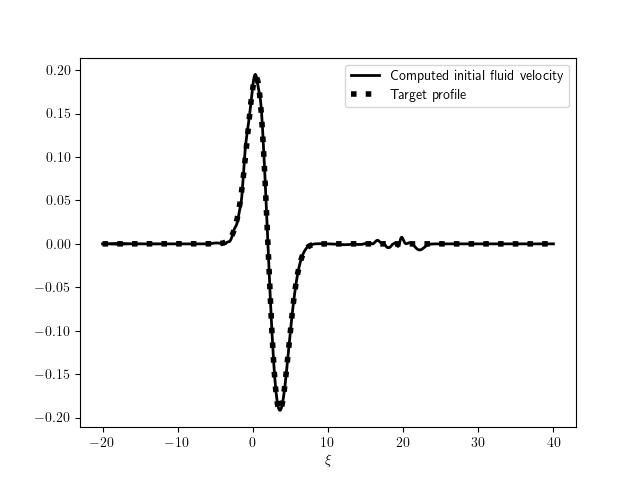}
		\caption{Iteration 2. }
	    \end{subfigure}
	\vfill
	     \begin{subfigure}{0.45\linewidth}
		 \includegraphics[width=\linewidth, height=\linewidth]{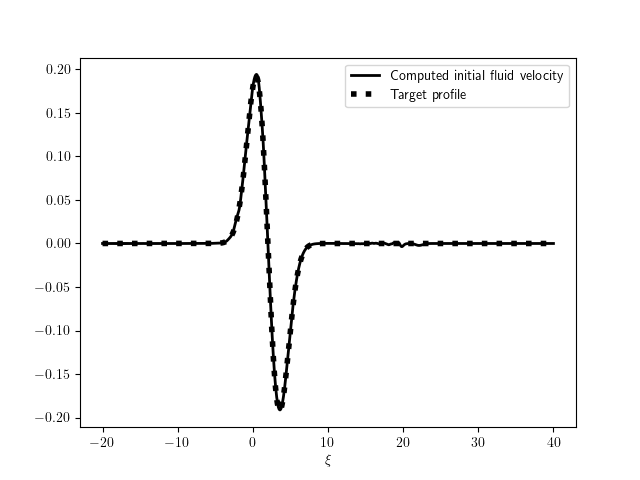}
		 \caption{Iteration 3. }
		 \label{}
	      \end{subfigure} 
	       \begin{subfigure}{0.45\linewidth}
		  \includegraphics[width=\linewidth, height=\linewidth ]{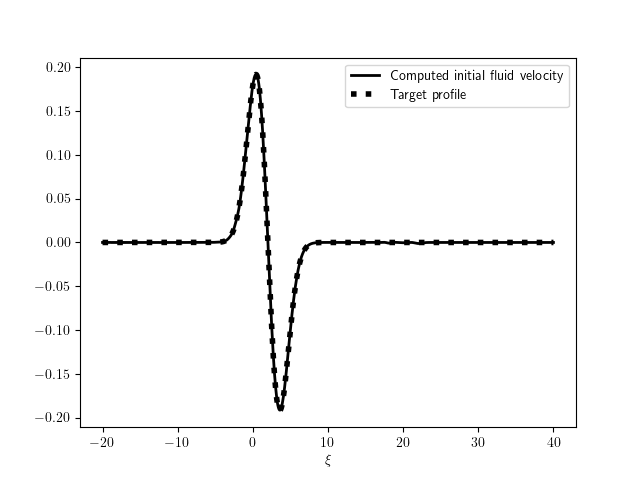}
		  \caption{Iteration 4. }
	       \end{subfigure}   
	\caption{Approximation of the initial fluid velocity $V_0(\xi)$ after the first four iterations of the L-BFGS-B algorithm used to minimize the functional \eqref{Functional_assimilation}.}
	\label{Iterations_V}
\end{figure}

Figures \ref{Iterations_N} and \ref{Iterations_V} show the results of the reconstruction procedure for both components $N$ and $V$, 
displaying the progression over the first four iterations of the L-BFGS-B scheme. The plots illustrate how the numerical solution evolves towards the target profiles $N_T$ and $V_T$, demonstrating the effectiveness of the variational assimilation strategy.

This experiment provides an initial numerical demonstration of the effectiveness of variational optimization for solving inverse problems involving the reconstruction of initial data in dispersive wave systems. Future research will focus on assessing the sensitivity of the reconstruction to observational noise, examining the influence of model parameters, and developing appropriate regularization strategies to improve the robustness and reliability of the method in practical scenarios.

\section{Conclusions}

In this work, we analyzed a nonlinear dispersive wave model with variable coefficients, capturing the dynamics of surface waves in a one-dimensional shallow water channel with uneven bottom topography. We established the well-posedness of the associated initial-value problem, ensuring that the model is mathematically robust and suitable for further analytical and numerical investigations.

To complement the theoretical results, we employed the numerical solver introduced in \cite{PizoMunoz}, which approximates solutions to the variable-coefficient Boussinesq system \eqref{Dee} over a finite spatial domain. Through a series of computational experiments, we investigated the system's response under varying parameter regimes and coefficient configurations, offering insight into the sensitivity and qualitative behavior of the model.

the numerical solver introduced in \cite{PizoMunoz}, which is capable of approximating solutions to the variable-coefficient Boussinesq system \eqref{Dee} over a finite spatial domain. Through a series of numerical experiments, we explored the system's behavior under different parameter regimes and coefficient variations, providing insight into the model's sensitivity and qualitative dynamics.

Furthermore, we initiated the study of an inverse problem associated with the system: the recovery of the initial wave profile and fluid velocity from final-time observations. This was approached by minimizing a discrepancy functional, offering a practical framework for estimating unknown initial conditions based on limited observational data.

These results lay the groundwork for a broader investigation of inverse problems in nonlinear dispersive systems. Future research may focus on the theoretical analysis of uniqueness and stability in the inverse setting, the extension to higher-dimensional domains, and the application of advanced optimization and machine learning techniques to enhance the reconstruction accuracy.

\section*{Acknowledgments}
This work was partially supported by Universidad del Valle, Cali, Colombia, under research projects C.I. 71235, 71288, Cali, Colombia, and MinCiencias under project FP44842-266-2017.

\end{document}